\documentclass[12pt]{amsart}

\usepackage{tikz}
\usetikzlibrary{calc, shapes.geometric, positioning, shapes, snakes, arrows, shapes, patterns}
\usepackage{caption}
\usepackage[labelformat=empty,labelfont={},margin=-.3em]{subcaption}

\usepackage{amsmath,amssymb,amsthm}
\usepackage{latexsym}
\usepackage[top=2in, bottom=1.5in, left=1in, right=1in]{geometry}
\usepackage{url}
\usepackage{color}

\newtheorem*{openproblem*}{Open Problem}

\newcommand{\newold}[2]{{\bf new:}{\color{blue}#1}{\bf old:}{\color{red}#2}}%
\newcommand{\w}{\mathsf{w}}
\renewcommand\thesubfigure{(\alph{subfigure})}
\newcommand\xxaxis{0}
\newcommand\yyaxis{90}
\newcommand\nsq[3]{\fill[fill=cyan!30, draw=black, rounded corners, line width=1pt, shift={(\xxaxis:#1)},shift={(\yyaxis:#2)}] (0,0) -- (2,0) -- (2,-2) -- (0,-2) --cycle;\node at (#1+1,#2-1) {{\scriptsize $#3$}};}

\begin{document}


\title[Commutation classes of the longest element in the symmetric group]{On the number of commutation classes of the longest element in the symmetric group} 

\author[Denoncourt]{Hugh Denoncourt}
\address{The City College of New York}
\email{hdenoncourt@ccny.cuny.edu}

\author[Ernst]{Dana C.~Ernst}
\author[Story]{Dustin Story}
\address{Northern Arizona University}
\email{dana.ernst@nau.edu, dustin@nau.edu}

\subjclass[2010]{05}


\begin{abstract}
Using the standard Coxeter presentation for the symmetric group $S_n$, two reduced expressions for the same group element are said to be commutation equivalent if we can obtain one expression from the other by applying a finite sequence of commutations. The resulting equivalence classes of reduced expressions are called commutation classes.  How many commutation classes are there for the longest element in $S_n$?

\bigskip

\emph{Original proposer of the open problem}: Donald E.~Knuth

\emph{The year when the open problem was proposed}: 1992~\cite[\S 9]{Knuth1992}

\end{abstract}


\maketitle

 
A \emph{Coxeter system} is a pair $(W,S)$ consisting of a distinguished (finite) set $S$ of generating involutions and a group
\[
W = \langle S \mid (st)^{m(s, t)} = e \text{ for } m(s, t) < \infty \rangle,
\]
called a \emph{Coxeter group}, where $e$ is the identity, $m(s,t) = 1$ if and only if $s = t$, and $m(s,t) = m(t,s)$. It turns out that the elements of $S$ are distinct as group elements and that $m(s, t)$ is the order of $st$.  Since the elements of $S$ have order two, the relation $(st)^{m(s,t)} = e$ can be written to allow the replacement 
\[
\underbrace{sts \cdots}_{m(s,t)} \mapsto \underbrace{tst \cdots}_{m(s,t)}
\]
which is called a \emph{commutation} if $m(s,t) = 2$ and a \emph{braid move} if $m(s,t) \geq 3$.

Given a Coxeter system $(W,S)$, a word $\w=s_{x_1}s_{x_2}\cdots s_{x_m}$ in the free monoid $S^*$ is called an \emph{expression} for $w\in W$ if it is equal to $w$ when considered as a group element. If $m$ is minimal among all expressions for $w$, the corresponding word is called a \emph{reduced expression} for $w$. In this case, we define the \emph{length} of $w$ to be $\ell(w)=m$. According to~\cite{Humphreys1990}, every finite Coxeter group contains a unique element of maximal length, which we refer to as the \emph{longest element} and denote by $w_0$.  

Let $(W,S)$ be a Coxeter system and let $w \in W$. Then $w$ may have several different reduced expressions that represent it. However, Matsumoto's Theorem~\cite[Theorem 1.2.2]{Geck2000} states that every reduced expression for $w$ can be obtained from any other by applying a finite sequence of commutations and braid moves.

Following~\cite{Stembridge1996}, we define a relation $\sim$ on the set of reduced expressions for $w$. Let $\w$ and $\w'$ be two reduced expressions for $w$ and define $\w \sim \w'$ if we can obtain $\w'$ from $\w$ by applying a single commutation. Now, define the equivalence relation $\approx$ by taking the reflexive transitive closure of $\sim$. Each equivalence class under $\approx$ is called a \emph{commutation class}.

The Coxeter system of type $A_{n-1}$ is generated by $S(A_{n-1}) = \{s_1, s_2, \ldots, s_{n-1}\}$ and has defining relations (i) $s_is_i = e$ for all $i$; (ii) $s_is_j = s_js_i$ when $|i-j| > 1$; and (iii) $s_is_js_i = s_js_is_j$ when $|i-j| = 1$. The corresponding Coxeter group $W(A_{n-1})$ is isomorphic to the symmetric group $S_{n}$ under the correspondence $s_i \mapsto (i,~i+1)$. 
%
It is well known that the longest element in $S_n$ is given in 1-line notation by
\[
w_0 = [n, n - 1, \ldots, 2, 1]
\]
and that $\ell(w_0)=\binom{n}{2}$.  

Let $c_n$ denote the number of commutation classes of the longest element in $S_n$.  The longest element $w_0$ in $S_4$ has length 6 and is given by the permutation $(1,4)(2,3)$.  There are 16 distinct reduced expressions for $w_0$ while $c_4=8$. The 8 commutation classes for $w_0$ are given in Figure~\ref{fig:CommClasses}, where we have listed the reduced expressions that each class contains.  Note that for brevity, we have written $i$ in place of $s_i$.

\begin{figure}
\begin{tabular}[b]{@{}cccccccc@{}}
\begin{tabular}{@{}c@{}}
321323\\
323123
\end{tabular} &
\begin{tabular}{@{}c@{}}
312312\\
132312\\
312132\\
132132
\end{tabular} &
\begin{tabular}{@{}c@{}}
321232
\end{tabular} &
\begin{tabular}{@{}c@{}}
232123
\end{tabular} &
\begin{tabular}{@{}c@{}}
123121\\
121321
\end{tabular} &
\begin{tabular}{@{}c@{}}
231231\\
213231\\
231213\\
213213
\end{tabular} &
\begin{tabular}{@{}c@{}}
123212
\end{tabular} &
\begin{tabular}{@{}c@{}}
212321
\end{tabular}
\end{tabular}
\caption{Reduced expressions and the corresponding commutation classes for the longest element in $S_4$.}\label{fig:CommClasses}
\end{figure}

In~\cite{Stanley1984}, Stanley provides a formula for the number of reduced expressions of the longest element $w_0$ in $S_n$.  However, the following question is currently unanswered.

\begin{openproblem*}
What is the number of commutation classes of the longest element in $S_n$?
\end{openproblem*}

To our knowledge, this problem was first introduced in 1992 by Knuth in Section 9 of~\cite{Knuth1992}, but not using our current terminology. A more general version of the problem appears in Section~5.2 of ~\cite{Kapranov1991}. In the paragraph following the proof of Proposition~4.4 of~\cite{Tenner2006}, Tenner explicitly states the open problem in terms of commutation classes. 

According to sequence A006245 of The On-Line Encyclopedia of Integer Sequences~\cite{OEIS}, the first 10 values for $c_n$ are $1, 1, 2, 8, 62, 908, 24698, 1232944, 112018190, 18410581880$.  To date, only the first 15 terms are known.  The current best upper-bound for $c_n$ was obtained by Felsner and Valtr.  They prove that for sufficiently large $n$, $c_n\leq 2^{0.6571n^2}$~\cite[Theorem 2]{Felsner2011}, although their result is stated in terms of arrangements of pseudolines.

The commutation classes of the longest element of the symmetric group are in bijection with a number of interesting objects.  It turns out that $c_n$ is equal to the number of 
\begin{itemize}
\item heaps for the longest element in $S_n$~\cite[Proposition 2.2]{Stembridge1996};
\item primitive sorting networks on $n$ elements~\cite{Armstrong2009,Kawahara2011,Knuth1992,Yamanaka2009,Yamanaka2010};\item rhombic tilings of a regular $2n$-gon (where all side lengths of the rhombi and the $2n$-gon are the same)~\cite{Elnitsky1997,Tenner2006};
\item oriented matroids of rank 3 on $n$ elements~\cite{Folkman1978,Kapranov1991};
\item arrangements of $n$ pseudolines~\cite{Felsner1997,Felsner2011,Knuth1992}.
\end{itemize}

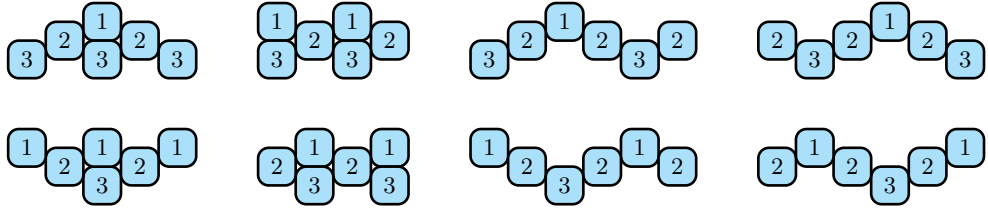
\begin{figure}
\centering
\begin{tabular}{llll}
\begin{tikzpicture}[scale=0.25]
\nsq{4}{3}{1}
\nsq{2}{2}{2}
\nsq{6}{2}{2}
\nsq{4}{1}{3}
\nsq{0}{1}{3}
\nsq{8}{1}{3}
\end{tikzpicture} \hspace{.3cm} &
\begin{tikzpicture}[scale=0.25]
\nsq{0}{3}{1}
\nsq{4}{3}{1}
\nsq{2}{2}{2}
\nsq{6}{2}{2}
\nsq{4}{1}{3}
\nsq{0}{1}{3}
\end{tikzpicture} \hspace{.3cm} &
\begin{tikzpicture}[scale=0.25]
\nsq{4}{3}{1}
\nsq{2}{2}{2}
\nsq{6}{2}{2}
\nsq{10}{2}{2}
\nsq{0}{1}{3}
\nsq{8}{1}{3}
\end{tikzpicture} \hspace{.3cm} &
\begin{tikzpicture}[scale=0.25]
\nsq{4}{3}{1}
\nsq{-2}{2}{2}
\nsq{2}{2}{2}
\nsq{6}{2}{2}
\nsq{0}{1}{3}
\nsq{8}{1}{3}
\end{tikzpicture} \\
\\
\begin{tikzpicture}[scale=0.25]
\nsq{4}{1}{3}
\nsq{2}{2}{2}
\nsq{6}{2}{2}
\nsq{4}{3}{1}
\nsq{0}{3}{1}
\nsq{8}{3}{1}
\end{tikzpicture} &
\begin{tikzpicture}[scale=0.25]
\nsq{0}{3}{1}
\nsq{4}{3}{1}
\nsq{2}{2}{2}
\nsq{-2}{2}{2}
\nsq{4}{1}{3}
\nsq{0}{1}{3}
\end{tikzpicture} &
\begin{tikzpicture}[scale=0.25]
\nsq{0}{3}{1}
\nsq{8}{3}{1}
\nsq{2}{2}{2}
\nsq{6}{2}{2}
\nsq{10}{2}{2}
\nsq{4}{1}{3}
\end{tikzpicture} &
\begin{tikzpicture}[scale=0.25]
\nsq{4}{1}{3}
\nsq{-2}{2}{2}
\nsq{2}{2}{2}
\nsq{6}{2}{2}
\nsq{0}{3}{1}
\nsq{8}{3}{1}
\end{tikzpicture}
\end{tabular}
\caption{Heaps for the longest element in $S_4$.}
\label{fig:Heaps}
\end{figure}

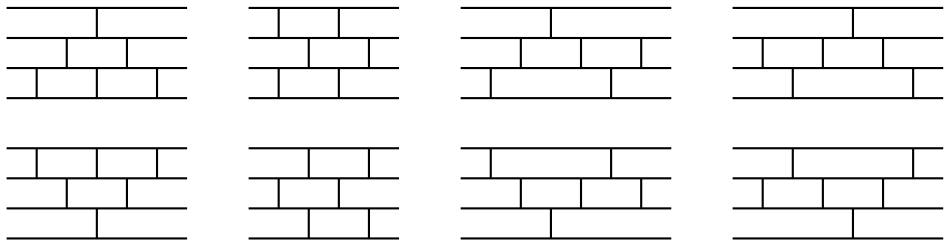
\begin{figure}
\centering
\begin{tabular}{llll}
\begin{tikzpicture}[scale=0.4]
\draw[thick] (0,0) -- (6,0);
\draw[thick] (0,1) -- (6,1);
\draw[thick] (0,2) -- (6,2);
\draw[thick] (0,3) -- (6,3);
\draw[thick] (1,0) -- (1,1);
\draw[thick] (3,0) -- (3,1);
\draw[thick] (5,0) -- (5,1);
\draw[thick] (2,1) -- (2,2);
\draw[thick] (4,1) -- (4,2);
\draw[thick] (3,2) -- (3,3);
\end{tikzpicture} \hspace{.3cm} &
\begin{tikzpicture}[scale=0.4]
\draw[thick] (0,0) -- (5,0);
\draw[thick] (0,1) -- (5,1);
\draw[thick] (0,2) -- (5,2);
\draw[thick] (0,3) -- (5,3);
\draw[thick] (1,0) -- (1,1);
\draw[thick] (3,0) -- (3,1);
\draw[thick] (1,2) -- (1,3);
\draw[thick] (2,1) -- (2,2);
\draw[thick] (4,1) -- (4,2);
\draw[thick] (3,2) -- (3,3);
\end{tikzpicture} \hspace{.3cm} &
\begin{tikzpicture}[scale=0.4]
\draw[thick] (0,0) -- (7,0);
\draw[thick] (0,1) -- (7,1);
\draw[thick] (0,2) -- (7,2);
\draw[thick] (0,3) -- (7,3);
\draw[thick] (1,0) -- (1,1);
\draw[thick] (6,1) -- (6,2);
\draw[thick] (5,0) -- (5,1);
\draw[thick] (2,1) -- (2,2);
\draw[thick] (4,1) -- (4,2);
\draw[thick] (3,2) -- (3,3);
\end{tikzpicture} \hspace{.3cm} &
\begin{tikzpicture}[scale=0.4]
\draw[thick] (0,0) -- (7,0);
\draw[thick] (0,1) -- (7,1);
\draw[thick] (0,2) -- (7,2);
\draw[thick] (0,3) -- (7,3);
\draw[thick] (2,0) -- (2,1);
\draw[thick] (1,1) -- (1,2);
\draw[thick] (6,0) -- (6,1);
\draw[thick] (3,1) -- (3,2);
\draw[thick] (5,1) -- (5,2);
\draw[thick] (4,2) -- (4,3);
\end{tikzpicture} \\
\\
\begin{tikzpicture}[scale=0.4]
\draw[thick] (0,0) -- (6,0);
\draw[thick] (0,1) -- (6,1);
\draw[thick] (0,2) -- (6,2);
\draw[thick] (0,3) -- (6,3);
\draw[thick] (1,2) -- (1,3);
\draw[thick] (3,2) -- (3,3);
\draw[thick] (5,2) -- (5,3);
\draw[thick] (2,1) -- (2,2);
\draw[thick] (4,1) -- (4,2);
\draw[thick] (3,0) -- (3,1);
\end{tikzpicture} &
\begin{tikzpicture}[scale=0.4]
\draw[thick] (0,0) -- (5,0);
\draw[thick] (0,1) -- (5,1);
\draw[thick] (0,2) -- (5,2);
\draw[thick] (0,3) -- (5,3);
\draw[thick] (2,0) -- (2,1);
\draw[thick] (4,0) -- (4,1);
\draw[thick] (2,2) -- (2,3);
\draw[thick] (3,1) -- (3,2);
\draw[thick] (1,1) -- (1,2);
\draw[thick] (4,2) -- (4,3);
\end{tikzpicture} &
\begin{tikzpicture}[scale=0.4]
\draw[thick] (0,0) -- (7,0);
\draw[thick] (0,1) -- (7,1);
\draw[thick] (0,2) -- (7,2);
\draw[thick] (0,3) -- (7,3);
\draw[thick] (1,2) -- (1,3);
\draw[thick] (6,1) -- (6,2);
\draw[thick] (5,2) -- (5,3);
\draw[thick] (2,1) -- (2,2);
\draw[thick] (4,1) -- (4,2);
\draw[thick] (3,0) -- (3,1);
\end{tikzpicture} &
\begin{tikzpicture}[scale=0.4]
\draw[thick] (0,0) -- (7,0);
\draw[thick] (0,1) -- (7,1);
\draw[thick] (0,2) -- (7,2);
\draw[thick] (0,3) -- (7,3);
\draw[thick] (2,2) -- (2,3);
\draw[thick] (1,1) -- (1,2);
\draw[thick] (6,2) -- (6,3);
\draw[thick] (3,1) -- (3,2);
\draw[thick] (5,1) -- (5,2);
\draw[thick] (4,0) -- (4,1);
\end{tikzpicture}
\end{tabular}
\caption{Minimal ladder lotteries corresponding to the primitive sorting networks on 4 elements.}
\label{fig:SortingNetworks}
\end{figure}

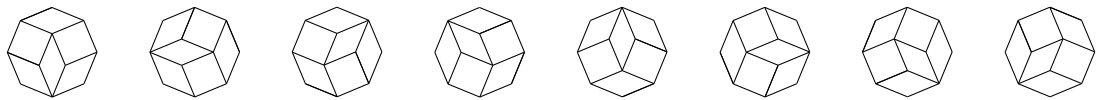
\begin{figure}
\centering
\begin{tabular}{llllllll}
\begin{tikzpicture}[scale=.25]
\draw [color=black] (1.0447355787000123,0.)-- (-0.6428174477092435,0.6928452343127731);
\draw [color=black] (-0.6428174477092435,0.6928452343127731)-- (-2.326013199790454,-0.010519390777748772);
\draw [color=black] (-2.326013199790454,-0.010519390777748772)-- (-3.018858434103227,-1.6980724171870047);
\draw [color=black] (-3.018858434103227,-1.6980724171870047)-- (-2.3154938090127053,-3.381268169268215);
\draw [color=black] (-2.3154938090127053,-3.381268169268215)-- (-0.6279407826034493,-4.074113403580988);
\draw [color=black] (-0.6279407826034493,-4.074113403580988)-- (1.0552549694777609,-3.3707487784904666);
\draw [color=black] (1.0552549694777609,-3.3707487784904666)-- (1.7481002037905344,-1.6831957520812106);
\draw [color=black] (1.7481002037905344,-1.6831957520812106)-- (1.0447355787000123,0.);
\draw (-2.326013199790454,-0.010519390777748772)-- (-0.6428174477092435,0.6928452343127731);
\draw (-2.326013199790454,-0.010519390777748772)-- (-0.6384601733811982,-0.7033646250905218);
\draw (-0.6384601733811982,-0.7033646250905218)-- (1.0447355787000123,0.);
\draw (-0.6384601733811982,-0.7033646250905218)-- (-1.3313054076939714,-2.390917651499778);
\draw (-3.018858434103227,-1.6980724171870047)-- (-1.3313054076939714,-2.390917651499778);
\draw (-1.3313054076939714,-2.390917651499778)-- (-0.6279407826034493,-4.074113403580988);
\draw (-0.6279407826034493,-4.074113403580988)-- (0.06490445170932387,-2.3865603771717327);
\draw (0.06490445170932387,-2.3865603771717327)-- (1.7481002037905344,-1.6831957520812106);
\draw (0.06490445170932387,-2.3865603771717327)-- (-0.6384601733811982,-0.7033646250905218);
\draw (-3.018858434103227,-1.6980724171870047)-- (-1.3313054076939714,-2.390917651499778);
\end{tikzpicture} \hspace{.2cm} &
\begin{tikzpicture}[scale=.25]
\begin{scope}[rotate=270]
\draw [color=black] (1.0447355787000123,0.)-- (-0.6428174477092435,0.6928452343127731);
\draw [color=black] (-0.6428174477092435,0.6928452343127731)-- (-2.326013199790454,-0.010519390777748772);
\draw [color=black] (-2.326013199790454,-0.010519390777748772)-- (-3.018858434103227,-1.6980724171870047);
\draw [color=black] (-3.018858434103227,-1.6980724171870047)-- (-2.3154938090127053,-3.381268169268215);
\draw [color=black] (-2.3154938090127053,-3.381268169268215)-- (-0.6279407826034493,-4.074113403580988);
\draw [color=black] (-0.6279407826034493,-4.074113403580988)-- (1.0552549694777609,-3.3707487784904666);
\draw [color=black] (1.0552549694777609,-3.3707487784904666)-- (1.7481002037905344,-1.6831957520812106);
\draw [color=black] (1.7481002037905344,-1.6831957520812106)-- (1.0447355787000123,0.);
\draw (-2.326013199790454,-0.010519390777748772)-- (-0.6428174477092435,0.6928452343127731);
\draw (-2.326013199790454,-0.010519390777748772)-- (-0.6384601733811982,-0.7033646250905218);
\draw (-0.6384601733811982,-0.7033646250905218)-- (1.0447355787000123,0.);
\draw (-0.6384601733811982,-0.7033646250905218)-- (-1.3313054076939714,-2.390917651499778);
\draw (-3.018858434103227,-1.6980724171870047)-- (-1.3313054076939714,-2.390917651499778);
\draw (-1.3313054076939714,-2.390917651499778)-- (-0.6279407826034493,-4.074113403580988);
\draw (-0.6279407826034493,-4.074113403580988)-- (0.06490445170932387,-2.3865603771717327);
\draw (0.06490445170932387,-2.3865603771717327)-- (1.7481002037905344,-1.6831957520812106);
\draw (0.06490445170932387,-2.3865603771717327)-- (-0.6384601733811982,-0.7033646250905218);
\draw (-3.018858434103227,-1.6980724171870047)-- (-1.3313054076939714,-2.390917651499778);
\end{scope}
\end{tikzpicture} \hspace{.2cm} &
\begin{tikzpicture}[scale=.25]
\begin{scope}[rotate=135]
\draw [color=black] (1.0447355787000123,0.)-- (-0.6428174477092435,0.6928452343127731);
\draw [color=black] (-0.6428174477092435,0.6928452343127731)-- (-2.326013199790454,-0.010519390777748772);
\draw [color=black] (-2.326013199790454,-0.010519390777748772)-- (-3.018858434103227,-1.6980724171870047);
\draw [color=black] (-3.018858434103227,-1.6980724171870047)-- (-2.3154938090127053,-3.381268169268215);
\draw [color=black] (-2.3154938090127053,-3.381268169268215)-- (-0.6279407826034493,-4.074113403580988);
\draw [color=black] (-0.6279407826034493,-4.074113403580988)-- (1.0552549694777609,-3.3707487784904666);
\draw [color=black] (1.0552549694777609,-3.3707487784904666)-- (1.7481002037905344,-1.6831957520812106);
\draw [color=black] (1.7481002037905344,-1.6831957520812106)-- (1.0447355787000123,0.);
\draw (-2.326013199790454,-0.010519390777748772)-- (-0.6428174477092435,0.6928452343127731);
\draw (-2.326013199790454,-0.010519390777748772)-- (-0.6384601733811982,-0.7033646250905218);
\draw (-0.6384601733811982,-0.7033646250905218)-- (1.0447355787000123,0.);
\draw (-0.6384601733811982,-0.7033646250905218)-- (-1.3313054076939714,-2.390917651499778);
\draw (-3.018858434103227,-1.6980724171870047)-- (-1.3313054076939714,-2.390917651499778);
\draw (-1.3313054076939714,-2.390917651499778)-- (-0.6279407826034493,-4.074113403580988);
\draw (-0.6279407826034493,-4.074113403580988)-- (0.06490445170932387,-2.3865603771717327);
\draw (0.06490445170932387,-2.3865603771717327)-- (1.7481002037905344,-1.6831957520812106);
\draw (0.06490445170932387,-2.3865603771717327)-- (-0.6384601733811982,-0.7033646250905218);
\draw (-3.018858434103227,-1.6980724171870047)-- (-1.3313054076939714,-2.390917651499778);
\end{scope}
\end{tikzpicture} \hspace{.2cm} &
\begin{tikzpicture}[scale=.25]
\begin{scope}[rotate=225]
\draw [color=black] (1.0447355787000123,0.)-- (-0.6428174477092435,0.6928452343127731);
\draw [color=black] (-0.6428174477092435,0.6928452343127731)-- (-2.326013199790454,-0.010519390777748772);
\draw [color=black] (-2.326013199790454,-0.010519390777748772)-- (-3.018858434103227,-1.6980724171870047);
\draw [color=black] (-3.018858434103227,-1.6980724171870047)-- (-2.3154938090127053,-3.381268169268215);
\draw [color=black] (-2.3154938090127053,-3.381268169268215)-- (-0.6279407826034493,-4.074113403580988);
\draw [color=black] (-0.6279407826034493,-4.074113403580988)-- (1.0552549694777609,-3.3707487784904666);
\draw [color=black] (1.0552549694777609,-3.3707487784904666)-- (1.7481002037905344,-1.6831957520812106);
\draw [color=black] (1.7481002037905344,-1.6831957520812106)-- (1.0447355787000123,0.);
\draw (-2.326013199790454,-0.010519390777748772)-- (-0.6428174477092435,0.6928452343127731);
\draw (-2.326013199790454,-0.010519390777748772)-- (-0.6384601733811982,-0.7033646250905218);
\draw (-0.6384601733811982,-0.7033646250905218)-- (1.0447355787000123,0.);
\draw (-0.6384601733811982,-0.7033646250905218)-- (-1.3313054076939714,-2.390917651499778);
\draw (-3.018858434103227,-1.6980724171870047)-- (-1.3313054076939714,-2.390917651499778);
\draw (-1.3313054076939714,-2.390917651499778)-- (-0.6279407826034493,-4.074113403580988);
\draw (-0.6279407826034493,-4.074113403580988)-- (0.06490445170932387,-2.3865603771717327);
\draw (0.06490445170932387,-2.3865603771717327)-- (1.7481002037905344,-1.6831957520812106);
\draw (0.06490445170932387,-2.3865603771717327)-- (-0.6384601733811982,-0.7033646250905218);
\draw (-3.018858434103227,-1.6980724171870047)-- (-1.3313054076939714,-2.390917651499778);
\end{scope}
\end{tikzpicture} \hspace{.2cm} &
\begin{tikzpicture}[scale=.25]
\begin{scope}[rotate=180]
\draw [color=black] (1.0447355787000123,0.)-- (-0.6428174477092435,0.6928452343127731);
\draw [color=black] (-0.6428174477092435,0.6928452343127731)-- (-2.326013199790454,-0.010519390777748772);
\draw [color=black] (-2.326013199790454,-0.010519390777748772)-- (-3.018858434103227,-1.6980724171870047);
\draw [color=black] (-3.018858434103227,-1.6980724171870047)-- (-2.3154938090127053,-3.381268169268215);
\draw [color=black] (-2.3154938090127053,-3.381268169268215)-- (-0.6279407826034493,-4.074113403580988);
\draw [color=black] (-0.6279407826034493,-4.074113403580988)-- (1.0552549694777609,-3.3707487784904666);
\draw [color=black] (1.0552549694777609,-3.3707487784904666)-- (1.7481002037905344,-1.6831957520812106);
\draw [color=black] (1.7481002037905344,-1.6831957520812106)-- (1.0447355787000123,0.);
\draw (-2.326013199790454,-0.010519390777748772)-- (-0.6428174477092435,0.6928452343127731);
\draw (-2.326013199790454,-0.010519390777748772)-- (-0.6384601733811982,-0.7033646250905218);
\draw (-0.6384601733811982,-0.7033646250905218)-- (1.0447355787000123,0.);
\draw (-0.6384601733811982,-0.7033646250905218)-- (-1.3313054076939714,-2.390917651499778);
\draw (-3.018858434103227,-1.6980724171870047)-- (-1.3313054076939714,-2.390917651499778);
\draw (-1.3313054076939714,-2.390917651499778)-- (-0.6279407826034493,-4.074113403580988);
\draw (-0.6279407826034493,-4.074113403580988)-- (0.06490445170932387,-2.3865603771717327);
\draw (0.06490445170932387,-2.3865603771717327)-- (1.7481002037905344,-1.6831957520812106);
\draw (0.06490445170932387,-2.3865603771717327)-- (-0.6384601733811982,-0.7033646250905218);
\draw (-3.018858434103227,-1.6980724171870047)-- (-1.3313054076939714,-2.390917651499778);
\end{scope}
\end{tikzpicture} \hspace{.2cm} &
\begin{tikzpicture}[scale=.25]
\begin{scope}[rotate=90]
\draw [color=black] (1.0447355787000123,0.)-- (-0.6428174477092435,0.6928452343127731);
\draw [color=black] (-0.6428174477092435,0.6928452343127731)-- (-2.326013199790454,-0.010519390777748772);
\draw [color=black] (-2.326013199790454,-0.010519390777748772)-- (-3.018858434103227,-1.6980724171870047);
\draw [color=black] (-3.018858434103227,-1.6980724171870047)-- (-2.3154938090127053,-3.381268169268215);
\draw [color=black] (-2.3154938090127053,-3.381268169268215)-- (-0.6279407826034493,-4.074113403580988);
\draw [color=black] (-0.6279407826034493,-4.074113403580988)-- (1.0552549694777609,-3.3707487784904666);
\draw [color=black] (1.0552549694777609,-3.3707487784904666)-- (1.7481002037905344,-1.6831957520812106);
\draw [color=black] (1.7481002037905344,-1.6831957520812106)-- (1.0447355787000123,0.);
\draw (-2.326013199790454,-0.010519390777748772)-- (-0.6428174477092435,0.6928452343127731);
\draw (-2.326013199790454,-0.010519390777748772)-- (-0.6384601733811982,-0.7033646250905218);
\draw (-0.6384601733811982,-0.7033646250905218)-- (1.0447355787000123,0.);
\draw (-0.6384601733811982,-0.7033646250905218)-- (-1.3313054076939714,-2.390917651499778);
\draw (-3.018858434103227,-1.6980724171870047)-- (-1.3313054076939714,-2.390917651499778);
\draw (-1.3313054076939714,-2.390917651499778)-- (-0.6279407826034493,-4.074113403580988);
\draw (-0.6279407826034493,-4.074113403580988)-- (0.06490445170932387,-2.3865603771717327);
\draw (0.06490445170932387,-2.3865603771717327)-- (1.7481002037905344,-1.6831957520812106);
\draw (0.06490445170932387,-2.3865603771717327)-- (-0.6384601733811982,-0.7033646250905218);
\draw (-3.018858434103227,-1.6980724171870047)-- (-1.3313054076939714,-2.390917651499778);
\end{scope}
\end{tikzpicture} \hspace{.2cm} &
\begin{tikzpicture}[scale=.25]
\begin{scope}[rotate=45]
\draw [color=black] (1.0447355787000123,0.)-- (-0.6428174477092435,0.6928452343127731);
\draw [color=black] (-0.6428174477092435,0.6928452343127731)-- (-2.326013199790454,-0.010519390777748772);
\draw [color=black] (-2.326013199790454,-0.010519390777748772)-- (-3.018858434103227,-1.6980724171870047);
\draw [color=black] (-3.018858434103227,-1.6980724171870047)-- (-2.3154938090127053,-3.381268169268215);
\draw [color=black] (-2.3154938090127053,-3.381268169268215)-- (-0.6279407826034493,-4.074113403580988);
\draw [color=black] (-0.6279407826034493,-4.074113403580988)-- (1.0552549694777609,-3.3707487784904666);
\draw [color=black] (1.0552549694777609,-3.3707487784904666)-- (1.7481002037905344,-1.6831957520812106);
\draw [color=black] (1.7481002037905344,-1.6831957520812106)-- (1.0447355787000123,0.);
\draw (-2.326013199790454,-0.010519390777748772)-- (-0.6428174477092435,0.6928452343127731);
\draw (-2.326013199790454,-0.010519390777748772)-- (-0.6384601733811982,-0.7033646250905218);
\draw (-0.6384601733811982,-0.7033646250905218)-- (1.0447355787000123,0.);
\draw (-0.6384601733811982,-0.7033646250905218)-- (-1.3313054076939714,-2.390917651499778);
\draw (-3.018858434103227,-1.6980724171870047)-- (-1.3313054076939714,-2.390917651499778);
\draw (-1.3313054076939714,-2.390917651499778)-- (-0.6279407826034493,-4.074113403580988);
\draw (-0.6279407826034493,-4.074113403580988)-- (0.06490445170932387,-2.3865603771717327);
\draw (0.06490445170932387,-2.3865603771717327)-- (1.7481002037905344,-1.6831957520812106);
\draw (0.06490445170932387,-2.3865603771717327)-- (-0.6384601733811982,-0.7033646250905218);
\draw (-3.018858434103227,-1.6980724171870047)-- (-1.3313054076939714,-2.390917651499778);
\end{scope}
\end{tikzpicture} \hspace{.2cm} &
\begin{tikzpicture}[scale=.25]
\begin{scope}[rotate=315]
\draw [color=black] (1.0447355787000123,0.)-- (-0.6428174477092435,0.6928452343127731);
\draw [color=black] (-0.6428174477092435,0.6928452343127731)-- (-2.326013199790454,-0.010519390777748772);
\draw [color=black] (-2.326013199790454,-0.010519390777748772)-- (-3.018858434103227,-1.6980724171870047);
\draw [color=black] (-3.018858434103227,-1.6980724171870047)-- (-2.3154938090127053,-3.381268169268215);
\draw [color=black] (-2.3154938090127053,-3.381268169268215)-- (-0.6279407826034493,-4.074113403580988);
\draw [color=black] (-0.6279407826034493,-4.074113403580988)-- (1.0552549694777609,-3.3707487784904666);
\draw [color=black] (1.0552549694777609,-3.3707487784904666)-- (1.7481002037905344,-1.6831957520812106);
\draw [color=black] (1.7481002037905344,-1.6831957520812106)-- (1.0447355787000123,0.);
\draw (-2.326013199790454,-0.010519390777748772)-- (-0.6428174477092435,0.6928452343127731);
\draw (-2.326013199790454,-0.010519390777748772)-- (-0.6384601733811982,-0.7033646250905218);
\draw (-0.6384601733811982,-0.7033646250905218)-- (1.0447355787000123,0.);
\draw (-0.6384601733811982,-0.7033646250905218)-- (-1.3313054076939714,-2.390917651499778);
\draw (-3.018858434103227,-1.6980724171870047)-- (-1.3313054076939714,-2.390917651499778);
\draw (-1.3313054076939714,-2.390917651499778)-- (-0.6279407826034493,-4.074113403580988);
\draw (-0.6279407826034493,-4.074113403580988)-- (0.06490445170932387,-2.3865603771717327);
\draw (0.06490445170932387,-2.3865603771717327)-- (1.7481002037905344,-1.6831957520812106);
\draw (0.06490445170932387,-2.3865603771717327)-- (-0.6384601733811982,-0.7033646250905218);
\draw (-3.018858434103227,-1.6980724171870047)-- (-1.3313054076939714,-2.390917651499778);
\end{scope}
\end{tikzpicture}
\end{tabular}
\caption{Rhombic tilings of a regular octagon.}
\label{fig:RhombicTilings}
\end{figure}

In Figure~\ref{fig:Heaps}, we have drawn lattice point representations of the 8 heaps that correspond to the commutation classes for the longest element in $S_4$.  Note that our heaps are sideways versions of the heaps that usually appear in the literature.  The minimum ladder lotteries (or ghost legs) corresponding to the 8 primitive sorting networks on $4$ elements are provided in Figure~\ref{fig:SortingNetworks}.  The 8 distinct rhombic tilings of a regular octagon are depicted in Figure~\ref{fig:RhombicTilings}.  

Very little is known about the number of commutation classes of the longest element in other finite Coxeter groups.


\bibliographystyle{plain}
\bibliography{CountingCommutationClasses}

\end{document}